\documentclass[12pt,twoside]{article}
\usepackage[english]{babel}
\usepackage[latin1]{inputenc}
\usepackage{amsmath}
\usepackage{amssymb,amsfonts}
\usepackage{graphicx}                   

\newcommand{\bG}{\mathbf{G}}

\newcommand{\bh}{\mathbf{h}}

\newcommand{\ba}{\mathbf{a}}
\newcommand{\bc}{\mathbf{c}}

\newcommand{\HYP}{\mathbb{H}^3}

\begin{document}
\pagestyle{myheadings}
\markboth{\centerline{Mikl\'os Eper and Jen\H o Szirmai}}
{Coverings with congruent and non-congruent hyperballs \dots}
\title
{Coverings with congruent and non-congruent hyperballs generated by doubly truncated Coxeter orthoschemes}

\author{\normalsize{Mikl\'os Eper and Jen\H o Szirmai} \\
\normalsize Budapest University of Technology and \\
\normalsize Economics Institute of Mathematics, \\
\normalsize Department of Geometry \\
\date{\normalsize{\today}}}

\maketitle
%
\begin{abstract}
After the investigation of the congruent and non-congruent hyperball packings related to 
doubly truncated Coxeter orthoscheme tilings \cite{SzJ1}, we consider the corresponding covering problems. 
In \cite{MSSz} the authors gave a partial classification of supergroups of some hyperbolic space groups whose fundamental domains will be integer parts of
truncated tetrahedra, and determined the optimal congruent hyperball packing and covering configurations 
belonging to some of these classes. 

In this paper we compliment these results with the investigation of the non-congruent covering cases, 
and the remainig congruent cases. We prove, that between congruent and non-congruent hyperball coverings the thinnest belongs to the $\{7,3,7\}$ Coxeter tiling with density 
$\approx 1.26829$. This covering density is smaller than the conjectured lower bound density of L.~Fejes~T\'oth for coverings with balls and horoballs.

We also study the local packing arrangements related to $\{u,3,7\}$ $(6< u < 7, ~ u\in \mathbb{R})$ doubly truncated orthoschemes and the corresponding hyperball coverings. 
We prove, that these coverings are achieved their minimum density at parameter $u\approx 6.45953$ with co- vering density $\approx 1.26454$ which is smaller then the above record-small
density, but this hyperball arrangement related to this locally optimal covering can not be extended to the entire $\mathbb{H}^3$.

Moreover we see, that in the hyperbolic plane $\mathbb{H}^2$ the universal lower bound of the congruent circle, horocycle, hypercycle covering density 
$\frac{\sqrt{12}}{\pi}$ can be approximated arbitrarily well also with non-congruent hypercycle coverings generated by doubly truncated Co- xeter orthoschemes.
\end{abstract}

\newtheorem{theorem}{Theorem}[section]
\newtheorem{corollary}[theorem]{Corollary}
\newtheorem{conjecture}{Conjecture}[section]
\newtheorem{lemma}[theorem]{Lemma}
\newtheorem{exmple}[theorem]{Example}
\newtheorem{defn}[theorem]{Definition}
\newtheorem{rmrk}[theorem]{Remark}
\newenvironment{definition}{\begin{defn}\normalfont}{\end{defn}}
\newenvironment{remark}{\begin{rmrk}\normalfont}{\end{rmrk}}
\newenvironment{example}{\begin{exmple}\normalfont}{\end{exmple}}
\newenvironment{acknowledgement}{Acknowledgement}


\section{Introduction}
The investigation of optimal dense packings and coverings with congruent balls in spaces with constant curvature, is one of the important topics in 
discrete geometry. The most explored is the Euclidean case, here one of the greatest results is the solution of the famous Kepler-conjecture 
\cite{Kep} ({\sc Hilbert}'s 18th problem). It's computer supported proof was given by {\sc Thomas Hales} in the early 2000s \cite{Hal}, 
which based on the ideas of {\sc L\' aszl\' o Fejes T\' oth} \cite{FTL}.

In the hyperbolic case there are far more open questions, in $n$-dimension $(n\geq 3)$ e.g. it isn't clear, when the most dense packing is realised with classical balls.
The so far known highest packing density is $\approx 0.77147$ with classical balls in $\mathbb{H}^3$, published in \cite{MSz}, where the authors also determined a 
classical ball configuration, which provides the so far known thinnest covering density $\approx 1.36893$.

Moreover, in hyperbolic space the definition of the desnsity of the packings and coverings is crucial, as it was shown by 
{\sc K\' aroly B\" or\" oczky} in his works \cite{BK1}, \cite{BK2}. The most common density-definition considers the local density of the balls related 
to their Dirichlet-Voronoi cells. We will use as well this local density-definition and also it's extension. Very important result of the classical 
ball and horoball packings the following theorem (see \cite{BK3}, \cite{BK4}):
\begin{theorem}[{\sc K.~B\"or\"oczky}]
In an $n$-dimensional space of constant curvature, consider a packing of spheres of radius $r$.
In the case of spherical space, assume that $r<\frac{\pi}{4}$.
Then the density of each sphere in its Dirichlet--Voronoi cell cannot exceed the density of $n+1$ spheres of radius $r$ mutually
touching one another with respect to the simplex spanned by their centers.
\end{theorem}
In hyperbolic space $\mathbb{H}^n$ $(n\geq 2)$ in addition to the classical spheres, there are two other types of balls: 
horoballs and hyperballs, which are non-compact "balls", and the above packing and covering problems with these kind of balls were also 
intensively investigated. The densest packing configuration, which Theorem 1.1 states, can be realized in $\mathbb{H}^3$, 
but surprisingly not with classical balls, but with horoballs, providing density $\approx 0.85328$ 
(this density can be attained with more different types of horoball packings related to fully assymptotic Coxeter tiling \cite{KSz}). 
In higher dimensions ($n=4,...,9$) there are also interesting results with high densities (see \cite{KSz2}, \cite{KSz3}, \cite{KSz4}). 
It was also shown, that in $\mathbb{H}^n$ ($n\geq 4$) the density of locally optimal horoball packing  related to $n$-dimensional fully assymptotic regular simplices exceeds 
the conjectured bound, for example in $\mathbb{H}^4$ it attains density $\approx 0.77147$, but the corresponding configuration 
can not be extended to the entire hyperbolic space (see \cite{SzJ2}, \cite{SzJ3}). 
Another question is, what will be the configuration in certain dimension for optimal horoball packing and covering with horoballs of 
"different types" \cite{KSz5}, \cite{SzJ4}. 

In the epoch-making book of {\sc L\' aszl\' o Fejes T\' oth} \cite{FTL} one can find the 
description of a horoball covering with density $\approx 1.280$, wich belongs to the $\{6,3,3\}$ tiling, and here conjectured that this would give the thinnest covering in the 
$\mathbb{H}^3$.

We know even less about the packings and coverings with hyperballs. 
In the hyperbolic plane ($n=2$) {\sc I. Vermes} proved, that the upper bound for {\it congruent hypercycle packing density} is $\frac{3}{\pi}$ \cite{VI1}, 
and the lower bound for {\it congruent hypercycle covering density} is $\frac{\sqrt{12}}{\pi}$ \cite{VI2}. We note, that the densest hypercycle 
packings and least dense covering can be realized 
with the same density as the optomal packing and coverings using horospheres, as it was shown e.g. in \cite{VI3}. 
However, there are no results related to the densest hypercycle packings and thinnest hypercycle coverings
with non-congruent hypercycles.

Moreover, in  higher dimensions there are few results related to congruent and non-congruent hyperball packings and coverings. 
The locally optimal hyperball packing configuration was researched in previous works related to several tilings: 
to tilings with truncated regular hyperbolic simplices in 
\cite{SzJ5}, \cite{SzJ6}; to tilings with cubes and octahedrons in \cite{SzJ7} 
(that yields density $\approx 0.84931$ in non-congruent cases); to tilings with truncated Coxeter orthoschemes 
in \cite{SzJ1} (that yields density $\approx 0.81335$). 

But these are just some of the works, where the densest hyperball packings were 
investigated in $\mathbb{H}^n$ (see also \cite{SzJ8}, \cite{SzJ9}). 
In \cite{SzJ10} the least dense hyperball 
coverings are determined related to 3-,4- and 5-dimensional Coxeter tilings, which can be derived 
from Coxeter orthoschemes by truncating the principal vertex with its polar plane. 
In several cases there were found locally optimal configurations, which densities exceeds the B\"or\"oczky-Florian 
upper bound related to classical ball and horoball packings. 

Another types of the ball packings and coverings in $\mathbb{H}^n$ are the so called 
hyp-hor packings \cite{SzJ11} and coverings \cite{ESz}, where in a configurations we use both horoball and hyperball, 
and the fundamental domain of the tilings are simply truncated simply asymptotic Coxeter orthoschemes. 
In the hyperbolic plane in both cases (packing and covering) the theoretic bound of the density (as above) can be attained in limit. 
In $\mathbb{H}^3$ the optimal packing density is $\approx 0.83267$, and the optimal covering density is $\approx 
1.27297$ both related to the $\{7,3,6\}$ Coxeter tiling. Moreover, we also considered configurations $\{p,3,6\}$ $(6< p < 7, ~ p\in \mathbb{R})$, 
and we got better densities: $\approx 0.85397$ for packing, $\approx 1.26885$ for covering, but these tilings can not be extended to the entire hyperbolic space.

In \cite{SzJ12} we considered congruent hyperball packings in
$3$-dimensional hyperbolic space and developed a decomposition algorithm that for each saturated hyperball packing provides a decomposition of $\HYP$
into truncated tetrahedra. Therefore, in order to get a density upper bound for hyperball packings, it is sufficient to determine
the density upper bound of hyperball packings in truncated simplices.
In \cite{SzJ13} we proved, using the above results, that the density upper bound of the saturated 
congruent hyperball (hypersphere) packings related to the
corresponding truncated tetrahedron cells is realized in regular truncated tetrahedra with density $\approx  0.86338$.
Furthermore, we prove that the density of locally optimal congruent hyperball arrangement 
in regular truncated tetrahedron is not monotonically increasing function of the height (radius) of corresponding optimal hyperball, 
contrary to the ball (sphere) and horoball (horosphere) packings.

Our discussion in this paper related to previous investigations (see \cite{MSSz}) where the authors considered some tilings generated by doubly truncated orthoschemes 
in $\mathbb{H}^3$, and determined their optimal packing and covering configurations with congruent hyperballs. 
In \cite{MSSz} was also shown, that the tiling $\{7,3,7\}$ provides 
$\approx 1.26829$ covering density, which is the currently known smallest ball-covering density in 
$\mathbb{H}^3$. We note here that in \cite{SzJ1} the correspon- ding packing problem was solved for congruent and non-congruent 
hyperballs related to similar orthoschemes. The tiling $\{7,3,7\}$ provides the optimal $\approx 0.81335$ packing density, that is realized with packings of congruent hyperballs.

In this paper, we complete and close the investigation started in the mentioned previous papers.
Now, we prove, that the thinnest covering with congruent or non-congruent hyperballs related to doubly truncated Coxeter orthoschemes generated tilings is realized 
at the $\{7,3,7\}$ tiling with density $\approx 1.26829$ (see also \cite{MSSz}).
{\it We note here that this covering density is smaller than the conjectured lower bound of L.~Fejes~T\'oth density for coverings with balls and horoballs.}

Moreover, we also consider $\{u,3,7\}$ $(6< u < 7, ~ u\in \mathbb{R})$ tilings, where the locally optimal 
density is $\approx 1.26454$, at parameter $u\approx 6.45953$. 

Finally we discuss the hypercycle coverings with congruent or non-congruent hypercycles, related to doubly truncated Coxeter orthoschemes generated tilings in 
the hyperbolic plane $\mathbb{H}^2$. We prove that the universal lower bound of the congruent circle or horocycle and hypercycle covering density 
$\frac{\sqrt{12}}{\pi}$ can be approximated arbitrarily well with non-congruent hypercycle coverings.
\section{Basic notions}
\subsection{The projective model of hyperbolic space $\mathbb{H}^3$}
For the computations we use the projective model \cite{ME} of the hyperbolic space. The model is defined in the $\mathbb{E}^{1,n}$ Lorentz space with signature $(1,n)$, i.e. consider $\mathbf{V}^{n+1}$ real vector space equipped with the bilinear form: 
\begin{equation*}
\langle ~ \mathbf{x},~\mathbf{y} \rangle = -x^0y^0+x^1y^1+ \dots + x^n y^n.
\end{equation*}
In the vector space consider the following equivalence relation:
\begin{equation*}
\mathbf{x}(x^0,...,x^n)\thicksim \mathbf{y}(y^0,...,y^n)\Leftrightarrow \exists c\in \mathbb{R}\backslash\{0\}: \mathbf{y}=c\cdot\mathbf{x}
\end{equation*}
The factorization with $\thicksim$ induces the $\mathcal{P}^n(\mathbf{V}^{n+1},\mbox{\boldmath$V$}\!_{n+1})$ $n$-dimensional real projective space. In this space to interpret the points of $\mathbb{H}^n$ hyperbolic space, consider the following quadratic form:
\begin{equation*}
Q=\{[\mathbf{x}]\in\mathcal{P}^n | \langle ~ \mathbf{x},~\mathbf{x} \rangle =0 \}=:\partial \mathbb{H}^n
\end{equation*}
The inner points relative to the cone-component determined by $Q$ are the points of $\mathbb{H}^n$ (for them $\langle ~ \mathbf{x},~\mathbf{x} \rangle <0$), the point of $Q=\partial \mathbb{H}^n$ are called the points at infinity, and the points lying outside realtive to $Q$ are outer points of $\mathbb{H}^n$ (for them $\langle ~ \mathbf{x},~\mathbf{x} \rangle >0$). We can also define a linear polarity between the points and hyperplanes of the space: the polar hyperplane of a point $[\mathbf{x}]\in\mathcal{P}^n$ is $Pol(\mathbf{x}):=\{[\mathbf{y}]\in\mathcal{P}^n | \langle ~ \mathbf{x},~\mathbf{y} \rangle =0 \}$, and hence $\mathbf{x}\in\mathbf{V}^{n+1}$ is incident with ${\boldsymbol{a}}\in\mbox{\boldmath$V$}\!_{n+1}$ iff $\langle ~ \mathbf{x},~{\boldsymbol{a}} \rangle =0$. In this projective model we can define a metric structure related to the above bilinear form, where for the distance of two proper points:
\begin{equation}
\cosh(d(\mathbf{x},\mathbf{y}))=\frac{-\langle ~ \mathbf{x},~\mathbf{y} \rangle}{\sqrt{\langle ~ \mathbf{x},~\mathbf{x} \rangle \langle ~ \mathbf{y},~\mathbf{y} \rangle}}.
\end{equation}
This corresponds to the distance formula in the well-known Beltrami-Cayley-Klein model.
\subsection{Hyperballs}
We implement the covering of $\mathbb{H}^n$ with hyperballs: we assign hyperballs to a doubly truncated orthoscheme, and if the hyperballs cover it, than accor- dingly to the images of the orthoscheme, the images of the hyperballs provide a covering of the space.

The equidistant surface (or hypersphere) is a quadratic surface that lies at a constant distance
from a plane in both halfspaces. The infinite body of the hypersphere is called a hyperball.
The $n$-dimensional {\it half-hypersphere } $(n=2,3)$ with distance $h$ to a hyperplane $\pi$
is denoted by $\mathcal{H}_n^h$. In the above model, the equation of $\mathcal{H}_n^h$ can be written \cite{SzP}:

\begin{equation}
1-\sum_{i=1}^{n}x_i^2=\left(-\frac{u_0}{\sinh h}+\sum_{i=1}^n \frac{u_i}{\sinh h}x_i\right)^2,
\end{equation}
where $(u_0,...,u_n)$ is the polar point of $\pi$ hyperplane.

The volume of a bounded hyperball piece $\mathcal{H}_n^h(\mathcal{A}_{n-1})$
bounded by an $(n-1)$-polytope $\mathcal{A}_{n-1} \subset \pi$, $\mathcal{H}_n^h$ and by
hyperplanes orthogonal to $\pi$ derived from the facets of $\mathcal{A}_{n-1}$ can be determined by the following formulas that follow from the suitable extension of the classical method of {\sc J. Bolyai} (\cite{BJ}):

\begin{equation}
Vol_3(\mathcal{H}_3^h(\mathcal{A}_2))=\frac{1}{4}Vol_2(\mathcal{A}_{2})\left[\sinh{(2h)}+
2 h \right],
\end{equation}
where the volume of the hyperbolic $(n-1)$-polytope $\mathcal{A}_{n-1}$ lying in the plane
$\pi$ is $Vol_{n-1}(\mathcal{A}_{n-1})$.

From the equation of the $\mathcal{H}_n^h$, we can see that the equivalents of the hyperballs in our model will be ellipsoids.
\section{Covering with hyperballs in hyperbolic space $\mathbb{H}^3$}
\subsection{Coxeter orthoschemes and tilings}

\begin{defn}
In the $\mathbb{H}^n$ $(2\leq n\in \mathbb{N})$ space a complete orthoscheme $\mathcal{O}$ of degree $d$ $(0\leq d\leq2)$ is a polytope bounded 
with hyperplanes $H^0,...,H^{n+d}$, for which $H^i\perp H^j$, unless $j\neq i-1,i,i+1$.
\end{defn}

In the classical ($d=0$) case let denote the vertex opposite to $H^i$ hyperplane with $A_i$ ($0\leq i\leq n$), and let denote the dihedral angle of 
$H^i$ and $H^j$ planes with $\alpha^{ij}$ (hence $\alpha^{ij}$ if $0\leq i<j-1\leq n$). 

In this paper we deal with orthoschemes of degree $d=2$, they can be described geometrically, as follows. 
We can give the sequence of the vertices of the orthoschemes $A_0,...,A_n$, where $A_iA_{i+1}$ edge is perpendicular 
to $A_{i+2}A_{i+3}$ edge for all $i\in\{0,...,n-3\}$. Here $A_0$ and $A_n$ are called the principal vertices of the orthoschemes. 
In the case $d=2$ these principal vertices are outer points of the above model, 
so they are truncated by its polar planes 
$Pol(A_0)$ and $Pol(A_n)$, and the orthoscheme is called doubly truncated (see Fig.~1). 
Now, we suppose that $A_0A_n$ line intersects the model, and don't deal with the other case.

In general the Coxeter orthoschemes were classified by {\sc H.-C. Im Hof}, he proved that they exist in dimension $\leq 9$, and gave a full list of them \cite{IH1} \cite{IH2}. 

\begin{figure}[ht]
\centering
\includegraphics[width=13.5cm]{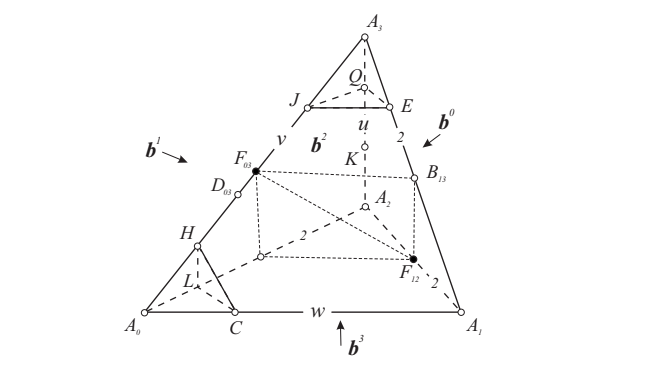}
\caption{3-dimensional Coxeter orthoscheme of degree $d=2$ with outer vertices $A_0,A_3$, truncated by $\pi ^0,\pi ^3 =HLC,JQE$ polar planes.}
\end{figure}

Now consider the reflections on the facets of the doubly truncated orthoscheme, and denote them with  $r_1,...,r_{n+3}$, hence define the group 
\begin{equation*}
\mathbf{G}=\langle r_1,...,r_{n+3} | (r_ir_j)^{m_{ij}}=1 \rangle,
\end{equation*}
where $\alpha^{ij}=\frac{\pi}{m_{ij}}$, so $m_{ii}=1$, and if $m_{ij}=\infty$ (i.e. $H^i$ and $H^j$ are parallel), than to the $r_i,r_j$ pair belongs no relation. 
Suppose that $2\leq m_{ij}\in\{\mathbb{N}\cup\infty\}$ if $i\neq j$. The Coxeter group $G$ acts on hyperbolic 
space $\overline{\mathbb{H}}^n$ properly discontinously, thus the images of the orthoscheme under this action provide a 
$\mathcal{T}$ tiling of $\overline{\mathbb{H}}^n$ (i.e. the images of the orthoscheme fills the $\overline{\mathbb{H}}^n$ without overlap).

For the {\it complete Coxeter orthoschemes} $\mathcal{O} \subset \mathbb{H}^n$ we adopt the usual
conventions and sometimes even use them in the Coxeter case: if two nodes are related by the weight $\cos{\frac{\pi}{m_{ij}}}$
then they are joined by a ($m_{ij}-2$)-fold line for $m_{ij}=3,~4$ and by a single line marked by $m_{ij}$ for $m_{ij} \geq 5$.
In the hyperbolic case if two bounding hyperplanes of $O$ are parallel, then the corresponding nodes
are joined by a line marked $\infty$. If they are divergent then their nodes are joined by a dotted line.

In the following we concentrate only on dimension $3$ and on hyperbolic
Coxeter-Schl\"afli symbol of the complete orthoscheme tiling $\mathcal{P}$ generated by reflections on the facets of a complete orthoscheme $\mathcal{O}$.
To every scheme there is a corresponding
symmetric $4 \times 4$ matrix $(b^{ij})$ where $b^{ii}=1$ and, for $i \ne j\in \{0,1,2,3\}$,
$b^{ij}$ equals to $-\cos{\alpha_{ij}}$ with all dihedral angles $\alpha_{ij}$ 
between the faces $H_i$,$H_j$ of $\mathcal{O}$.

For example, $(b^{ij})$ in formula (4) is the so called Coxeter-Schl\"afli matrix with
parameters $(u;v;w)$, i.e. $\alpha_{01}=\frac{\pi}{u}$, 
$\alpha_{12}=\frac{\pi}{v}$, $\alpha_{23}=\frac{\pi}{w}$.
Now only $3\le u,v,w$ come into account (see \cite{IH1,IH2}).
\begin{equation}
(b^{ij})=\langle \mbox{\boldmath$b^i$},\mbox{\boldmath$b^j$} \rangle:=\begin{pmatrix}
1& -\cos{\frac{\pi}{u}}& 0 & 0 \\
-\cos{\frac{\pi}{u}} & 1 & -\cos{\frac{\pi}{v}}& 0 \\
0 & -\cos{\frac{\pi}{v}} & 1 & -\cos{\frac{\pi}{w}} \\
0 & 0 & -\cos{\frac{\pi}{w}} & 1
\end{pmatrix}. 
\end{equation}

This $3$-dimensional complete (truncated or frustum) orthoscheme $\mathcal{O}=\mathcal{O}(u,v,w)$ and its reflection group $\bG_{uvw}$ will be described in
Fig.~1,~2, and by the symmetric Coxeter-Schl\"afli matrix $(b^{ij})$ in formula (4), furthermore by its inverse matrix $(h_{ij})$ in formula (5).
\begin{equation}
\begin{gathered}
(h_{ij})=(b^{ij})^{-1}=\langle \ba_i, \ba_j \rangle:=\\
=\frac{1}{B} \begin{pmatrix}
\sin^2{\frac{\pi}{w}}-\cos^2{\frac{\pi}{v}}& \cos{\frac{\pi}{u}}\sin^2{\frac{\pi}{w}}& \cos{\frac{\pi}{u}}\cos{\frac{\pi}{v}} & \cos{\frac{\pi}{u}}\cos{\frac{\pi}{v}}\cos{\frac{\pi}{w}} \\
\cos{\frac{\pi}{u}}\sin^2{\frac{\pi}{w}} & \sin^2{\frac{\pi}{w}} & \cos{\frac{\pi}{v}}& \cos{\frac{\pi}{w}}\cos{\frac{\pi}{v}} \\
\cos{\frac{\pi}{u}}\cos{\frac{\pi}{v}} & \cos{\frac{\pi}{v}} & \sin^2{\frac{\pi}{u}}  & \cos{\frac{\pi}{w}}\sin^2{\frac{\pi}{u}}  \\
\cos{\frac{\pi}{u}}\cos{\frac{\pi}{v}}\cos{\frac{\pi}{w}}  & \cos{\frac{\pi}{w}}\cos{\frac{\pi}{v}} & \cos{\frac{\pi}{w}}\sin^2{\frac{\pi}{u}}  & \sin^2{\frac{\pi}{u}}-\cos^2{\frac{\pi}{v}}
\end{pmatrix}, 
\end{gathered}
\end{equation}
where
$$
B=\det(b^{ij})=\sin^2{\frac{\pi}{u}}\sin^2{\frac{\pi}{w}}-\cos^2{\frac{\pi}{v}} <0, \ \ \text{i.e.} \ \sin{\frac{\pi}{u}}\sin{\frac{\pi}{w}}-\cos{\frac{\pi}{v}}<0.
$$
The volume of an doubly truncated $HLCA_1A_2JQE$ Coxeter orthoscheme with outer vertices $A_0,A_3$ (see Fig.~1) in $\mathbb{H}^n$ can be determined by the 
following theorem of {\sc R. Kellerhals} \cite{KH1,KH2}.
\begin{theorem} The volume of a three-dimensional hyperbolic
complete ortho\-scheme (except Lambert cube cases) $\mathcal{S}$
is expressed with the essential angles $\alpha^{01},\alpha^{12},\alpha^{23}, \ (0 \le \alpha^{ij} \le \frac{\pi}{2})$ in the following form:

\begin{align*}
&Vol_3(\mathcal{S})=\frac{1}{4} \{ \mathcal{L}(\alpha^{01}+\theta)-
\mathcal{L}(\alpha^{01}-\theta)+\mathcal{L}(\frac{\pi}{2}+\alpha^{12}-\theta)+ \notag \\
&+\mathcal{L}(\frac{\pi}{2}-\alpha^{12}-\theta)+\mathcal{L}(\alpha^{23}+\theta)-
\mathcal{L}(\alpha^{23}-\theta)+2\mathcal{L}(\frac{\pi}{2}-\theta) \}, 
\end{align*}
where $\theta \in [0,\frac{\pi}{2})$ is defined by the following formula:
$$
\tan(\theta)=\frac{\sqrt{ \cos^2{\alpha^{12}}-\sin^2{\alpha^{01}} \sin^2{\alpha^{23}
}}} {\cos{\alpha^{01}}\cos{\alpha^{23}}}
$$
and where $\mathcal{L}(x):=-\int\limits_0^x \log \vert {2\sin{t}} \vert dt$ \ denotes the
Lobachevsky function.
\end{theorem}
(Don't interest the Lambert cube case, when $A_0A_n$ goes out of the model.)
\subsection{Hyperball coverings}
Consider in 3-dimensional hyperbolic space $\mathbb{H}^3$ a doubly truncated Coxeter orthoscheme, with Coxeter-Schl\" afli symbol $\{u,v,w\}$. 
Its Coxeter-Schl\" afli matrix ($c^{ij}$) and its inverse are described in (4) and (5). 

Set the $A_0A_1A_2A_3$ orthoscheme in the model centered at $O=(1,0,0,0)$, so the $A_k[\mathbf{a_k}]$ $(k=1,2)$ vertices are proper points, 
i.e. $h_{kk}=\langle\mathbf{a}_k,\mathbf{a}_k\rangle<0$ if ($k=1,2$), and we can check easily that it's always fulfilled. 
On the other hand $A_k[\mathbf{a_k}]$ $(k=0,3)$ principal vertices are outer points, i.e. $h_{kk}=\langle\mathbf{a}_k,\mathbf{a}_k\rangle>0$ if ($k=0,3$), 
and from this we get the conditions $\sin^2{\frac{\pi}{w}}-\cos^2{\frac{\pi}{v}}<0$ and $\sin^2{\frac{\pi}{u}}-\cos^2{\frac{\pi}{v}}<0$, or equivalently 
$\frac{1}{w}+\frac{1}{v}<\frac{1}{2}$ and $\frac{1}{u}+\frac{1}{v}<\frac{1}{2}$. It means, that we have the following infinite series of 
$\mathcal{F}_u^{(v,w)}$ Coxeter orthoschemes with two outer vertices (see \cite{MSz} for details):

\begin{itemize}
\item[-]$\{u,v,w\}$, where $u\geq 3$, $v\geq 7$ and $w\geq 3$
\item[-]$\{u,v,w\}$, where $u\geq 4$, $v=5,6$ and $w\geq 4$
\item[-]$\{u,v,w\}$, where $u\geq 5$, $v=4$ and $w\geq 5$
\item[-]$\{u,v,w\}$, where $u\geq 7$, $v=3$ and $w\geq 7$
\end{itemize}

We truncate the orthoscheme with the polar planes $\pi^0= pol(A_0)[{\boldsymbol{a}^0}]$ and $\pi^3= pol(A_3)[{\boldsymbol{a}^3}]$ 
of the verties $A_0$ and $A_3$  (see Fig.~1), so $J=\pi^3\cap A_0A_3,Q=\pi^3\cap A_2A_3,E=\pi^3\cap A_1A_3$, and $H=\pi^0\cap A_0A_3,L=\pi^0\cap A_0A_2,C=\pi^0\cap A_0A_1$ 
are proper points. $Q$ lies on edge $A_2A_3$ so we can write ${\bf{q}}\sim c\cdot{\bf{a}}_3+{\bf{a}}_2$ for some $c$ real number. 
The corresponding $Q$ point lies on ${\boldsymbol{a}^3}$ iff their scalar product is $0$:
\begin{equation*}
c\cdot{\bf{a}}_3{\boldsymbol{a}^3}+{\bf{a}}_2{\boldsymbol{a}^3}=0 \Leftrightarrow c=-\frac{{\bf{a}}_2{\boldsymbol{a}^3}}{{\bf{a}}_3{\boldsymbol{a}^3}}
\end{equation*}
\begin{equation}
{\bf{q}}\sim -\frac{{\bf{a}}_2{\boldsymbol{a}^3}}{{\bf{a}}_3{\boldsymbol{a}^3}}{\bf{a}}_3+{\bf{a}}_2\sim {\bf{a}}_2({\bf{a}}_3{\boldsymbol{a}^3})-{\bf{a}}_3({\bf{a}}_2{\boldsymbol{a}^3})={\bf{a}}_2 h_{33}-{\bf{a}}_3 h_{23}
\end{equation}
Similarly for the other vertices on $\pi^3$ and $\pi^0$ polar planes.
\begin{equation}
\begin{gathered}
{\bf{j}}\sim {\bf{a}}_0({\bf{a}}_3{\boldsymbol{a}^3})-{\bf{a}}_3({\bf{a}}_0{\boldsymbol{a}^3})={\bf{a}}_0 h_{33}-{\bf{a}}_3 h_{03}
\\ {\bf{e}}\sim {\bf{a}}_1({\bf{a}}_3{\boldsymbol{a}^3})-{\bf{a}}_3({\bf{a}}_1{\boldsymbol{a}^3})={\bf{a}}_1 h_{33}-{\bf{a}}_3 h_{13}
\\{\bf{h}}\sim {\bf{a}}_3({\bf{a}}_0{\boldsymbol{a}^0})-{\bf{a}}_0({\bf{a}}_3{\boldsymbol{a}^0})={\bf{a}}_3 h_{0}-{\bf{a}}_0 h_{03}
\\ {\bf{l}}\sim {\bf{a}}_2({\bf{a}}_0{\boldsymbol{a}^0})-{\bf{a}}_0({\bf{a}}_2{\boldsymbol{a}^0})={\bf{a}}_2 h_{00}-{\bf{a}}_0 h_{02}
\\ {\bf{c}}\sim {\bf{a}}_1({\bf{a}}_0{\boldsymbol{a}^0})-{\bf{a}}_0({\bf{a}}_1{\boldsymbol{a}^0})={\bf{a}}_1 h_{0}-{\bf{a}}_0 h_{01}
\end{gathered} 
\end{equation}
Set the orthoscheme in the model with coordinates $Q=(1,0,0,0)$, $E=(1,0,y,0)$, $J=(1,x,y,0)$, $A_0=(1,x,y,-z_0)$, $A_1=(1,0,y,-z_1)$, $A_2=(1,0,0,-z_2)$, $H=(1,x,y,-z_H)$, 
$L=(t_1x,t_1y,-t_1z_0-(1-t_1)z_2)$, $C=(t_2x,y,-t_2z_0-(1-t_2)z_1)$, for some $t_1,t_2\in[0,1]$ (see Fig.~2), and using formulas (1) and (6-7):

\begin{gather*}
\cosh(d(Q,E))=\frac{-\left<{\bf{q}},{\bf{e}}\right>}{\sqrt{\left<{\bf{q}},{\bf{q}}\right>\left<{\bf{e}},{\bf{e}}\right>}}=\frac{h_{13}h_{23}-h_{12}h_{33}}{\sqrt{(h_{22}h_{33}-h_{23}^2)(h_{11}h_{33}-h_{13}^2)}}
\end{gather*}
\begin{gather*}
\cosh(d(Q,J))=\frac{-\left<{\bf{q}},{\bf{j}}\right>}{\sqrt{\left<{\bf{q}},{\bf{q}}\right>\left<{\bf{j}},{\bf{j}}\right>}}=\frac{h_{03}h_{23}-h_{02}h_{33}}{\sqrt{(h_{22}h_{33}-h_{23}^2)(h_{0}h_{33}-h_{03}^2)}}
\end{gather*}
\begin{gather*}
\cosh(d(E,A_1))=\frac{-\left<{\bf{a}}_1,{\bf{e}}\right>}{\sqrt{\left<{\bf{e}},{\bf{e}}\right>\left<{\bf{a}}_1,{\bf{a}}_1\right>}}=\sqrt{\frac{h_{11}h_{33}-h_{13}^2}{h_{11}h_{33}}}
\end{gather*}
\begin{gather*}
\cosh(d(Q,A_2))=\frac{-\left<{\bf{a}}_2,{\bf{q}}\right>}{\sqrt{\left<{\bf{q}},{\bf{q}}\right>\left<{\bf{a}}_2,{\bf{a}}_2\right>}}=\sqrt{\frac{h_{22}h_{33}-h_{23}^2}{h_{22}h_{33}}}
\end{gather*}
\begin{gather*}
\cosh(d(J,H))=\frac{-\left<{\bf{j}},{\bf{h}}\right>}{\sqrt{\left<{\bf{j}},{\bf{j}}\right>\left<{\bf{h}},{\bf{h}}\right>}}=\frac{-h_{02}h_{03}h_{33}+h_{03}^2h_{23}}{\sqrt{h_{00}h_{33}(h_{22}h_{33}-h_{23}^2)(h_{00}h_{33}-h_{03}^2)}}
\end{gather*}

We know furthermore, that $0=\left<\bh,\ba_0\right>=\left<\bc,\ba_0\right>=\left<\mathbf{l},\ba_0\right>$, hence we 
can determine the coordinates of the vertices of a certain $\mathcal{F}_u^{(v,w)}$.

\begin{figure}[ht]
\centering
\includegraphics[width=12.0cm]{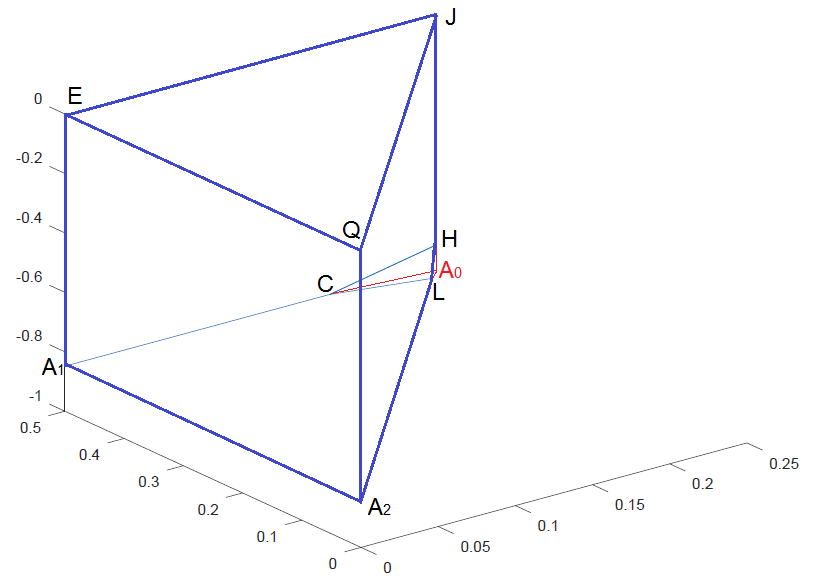}
\caption{$\{7,3,7\}$ Coxeter orthoscheme: $A_0$ vertex truncated with $HLC$, $A_3$ vertex truncated with $QEJ$ polar plane.}
\end{figure}

The images of this $\mathcal{F}_u^{(v,w)}$ doubly truncated orthoscheme under reflections on its facets provide the Coxeter tiling $\mathcal{T}_u^{(v,w)}$ 
with fundamental domain $\mathcal{F}_u^{(v,w)}$. We construct hyperball coverings as follows (see Fig.~4):

\begin{itemize}
\item[-] Let $QEJ$ the base hyperplane of a hyperball, and consider its piece bounded with $QEJ$ plane, 
the planes perpendicular to this base hyperplane, and contain the edges $QE,QJ,JE$. Denote this hyperball-piece with $\mathcal{H}_1$, and its height parameter with $h_1$.
\item[-] Let $HLC$ the base hyperplane of the other hyperball, and consider its piece bounded with $HLC$ plane, the planes perpendicular to this base hyperplane, 
and contain the edges $HL,HC,CL$. Denote this hyperball-piece with $\mathcal{H}_2$, and its height parameter with $h_2$.
\end{itemize}

It's obvious, that this hyperballs cover the orthoscheme iff they cover all of its edges, that is what we will check 
in the different cases. So if the hyperballs cover the edges of $\mathcal{F}_u^{(v,w)}$, the images of $\mathcal{H}_1$ and $\mathcal{H}_1$ 
under the reflections on the orthoschemes facets will provide a $\mathcal{C}_u^{(v,w)}$ covering of $\mathbb{H}^3$.

Here we remark, that the density can't be optimal, if we use only one hyperball to cover the orthoscheme, i.e. $h_i=0$ for the other one. 

A necessary condition for optimality is that there is a certain point of one of the $QA_2,EA_1,JH,LA_2,CA_1,A_1A_2$ edges 
lying on surfaces of both hyperballs. 
Therefore, the optimal covering density will be realized if the above conditions stand. According to this, we distinguish 6 cases.
\begin{defn}
The density of $\mathcal{C}_u^{(v,w)}$ covering:
\begin{equation*}
\delta(\mathcal{C}_u^{(v,w)})=\frac{Vol(\mathcal{H}_1)+Vol(\mathcal{H}_2)}{Vol(\mathcal{F}_u^{(v,w)})}
\end{equation*}
\end{defn}
\subsection{To cover, or not to cover?}
In order to decide whether a covering can be realized or not, we define the distance functions of the edges related to polar hyperplanes $\pi^0= pol(A_0)[{\boldsymbol{a}^0}]$ and $\pi^3= pol(A_3)[{\boldsymbol{a}^3}]$ of the orthoscheme. 
For example consider $HLC$ plane, $JH$ edge (see Fig.~2), and parametrize $JH$ by the following way: 
let $T$ a general point of $JH$, where $T(1,x,y,-tz_H),~t\in[0,1]$. Hence the distance function of $JH$ from $HLC$ at $t$, is the distance of $T$ and $H$, so by (1):
\begin{equation*}
d_{HLC}^{JH}(t)=\arccos\left(\frac{1-x^2-y^2-tz_H^2}{\sqrt{(1-x^2-y^2-z_H^2)(1-x^2-y^2-t^2z_H^2)}}\right)
\end{equation*}
Similarly we get functions $d_{HLC}^{LA_2}(t),d_{HLC}^{CA_1}(t),d_{HLC}^{QA_2}(t),d_{HLC}^{EA_1}(t),d_{HLC}^{A_1A_2}(t)$. In the other case consider $QEJ$ plane, and $JH$ edge for example. Here the function is given by the distance of $T(1,x,y,-tz_H),~t\in[0,1]$, and $J$:
\begin{equation*}
d_{QEJ}^{JH}(t)=\arccos\left(\sqrt{\frac{{1-x^2-y^2}}{{1-x^2-y^2-(tz_H)^2}}}\right)
\end{equation*}
And similarly for functions $d_{QEJ}^{QA_2}(t), d_{QEJ}^{EA_1}(t), d_{QEJ}^{LA_2}(t), d_{QEJ}^{CA_1}(t), d_{QEJ}^{A_1A_2}(t)$.

As we saw at the end of the previous section, we get optimal covering density, only if the surfaces of the two hyperballs intersect each other on an edge. 
So we choose an edge $e$ (which is not on one of the polar planes), 
and one of its point $T(t)$ (parametrized with $t\in [0,1]$), and say, that the $T$ lies on the surface of both hyperballs. 
Hence we know the $h_1,h_2$ heights of the hyperballs, and we have to check, whether the distance of the points of the other edges from one of the hyperplanes are 
smaller than the corresponding $h_i$ or not. 
We can determine the intersection points of a hyperball and an edge by using (2), and solving equations. 
If all of the points of an edge are closer to $QEJ$ than $h_1$ or to $HLC$ than $h_2$, than the covering is realized at $T$, 
and if this stands for all $T$ in $e$, than the covering is realized at edge $e$ ($e\in \{QA_2,EA_1,JH,LA_2,CA_1,A_1A_2\}$).

By the careful analysis of the above distance functions, using the help of computerial computations, we can say the followings:
\begin{itemize}
\item[-] The covering is realized at $A_1A_2$ edge (see Fig.~4).
\item[-] The covering is realized at $QA_2$ edge.
\item[-] The covering is realized at $CA_1$ edge.
\item[-] The covering isn't realized at $EA_1$ edge (the hyperballs don't cover $QA_2$ or $A_1A_2$ edge).
\item[-] The covering isn't realized at $LA_2$ edge (the hyperballs don't cover $CA_1$ or $A_1A_2$ edge).
\item[-] The covering isn't realized at $HJ$ edge.
\end{itemize}
\subsection{Non-congruent coverings}
In this subsection, we consider the non-congruent coverings for the possible cases described in the previous subsection.

We can determine the areas of $QEJ$ and $HLC$ triangulars, and if we settle $T(t)$ on the edge $QA_2$, than the heights of the hyperballs are $h_1(t)=d_{QEJ}^{QA_2}(t)$ and 
$h_2(t)=d_{HLC}^{QA_2}(t)$, thus we can compute the density of the covering using (3), Theorem 3.2, and Definition 3.3. 
The covering density will be a function with variable $t$, and we can determine its minimum precisely by real analysis (see Fig.~3). 
The following table shows the minimal covering densities at $QA_2$ for different types of orthoschemes. 

There are infinitely many types of doubly truncated orthoschemes, as we saw previously. Here we listed only some of the first elements of the four infinite 
sequences described in Section 3. With the further increase of the parameters $u,v,w$, the density grows as we obtain it after careful analysis of the density function.

\begin{center}
\begin{tabular}{|c|c|c|c|}
\hline
\multicolumn{1}{|c}{Type of orthoscheme} & \multicolumn{1}{|c|}{$\delta_{min}$} & \multicolumn{1}{|c|}{$h_1$} & \multicolumn{1}{|c|}{$h_2$}\\
\hline
$\mathcal{F}_3^{(7,3)}$      &$1.28943$ &$0.92295$ &$1.55521$ \\
\hline
$\mathcal{F}_3^{(8,3)}$      &$1.34248$ &$0.67445$ &$1.35737$ \\
\hline
$\mathcal{F}_4^{(5,4)}$      &$1.54311$ &$0.73337$ &$1.51710$ \\
\hline
$\mathcal{F}_4^{(6,4)}$      &$1.66605$ &$0.52867$ &$1.37017$ \\
\hline
$\mathcal{F}_5^{(4,5)}$      &$1.79576$ &$0.77124$ &$1.66724$ \\
\hline
$\mathcal{F}_5^{(5,4)}$      &$2.00292$ &$0.42347$ &$1.79770$ \\
\hline
$\mathcal{F}_6^{(4,5)}$      &$2.23585$ &$0.53126$ &$1.87500$ \\
\hline
$\mathcal{F}_6^{(5,4)}$      &$2.60090$ &$0.31440$ &$2.00574$ \\
\hline
$\mathcal{F}_7^{(3,7)}$      &$2.31671$ &$1.08534$ &$2.14790$ \\
\hline
$\mathcal{F}_7^{(4,5)}$      &$2.77700$ &$0.42041$ &$2.04284$ \\
\hline
\end{tabular}
\end{center}
The case of covering at $CA_1$ is the same, because of the symmetry of the orthoscheme. And see the optimal densities of the covering at the edge $A_1A_2$ in the following table.

\begin{center}
\begin{tabular}{|c|c|c|c|}
\hline
\multicolumn{1}{|c}{Type of orthoscheme} & \multicolumn{1}{|c|}{$\delta_{min}$} & \multicolumn{1}{|c|}{$h_1$} & \multicolumn{1}{|c|}{$h_2$}\\
\hline
$\mathcal{F}_3^{(7,3)}$      &$1.38712$ &$1.36405$ &$1.36405$ \\
\hline
$\mathcal{F}_3^{(8,3)}$      &$1.45345$ &$1.15039$ &$1.15039$ \\
\hline
$\mathcal{F}_4^{(5,4)}$      &$1.36411$ &$1.16974$ &$1.16974$ \\
\hline
$\mathcal{F}_4^{(5,5)}$      &$1.41055$ &$1.29237$ &$0.85103$ \\
\hline
$\mathcal{F}_5^{(4,5)}$      &$1.31751$ &$1.19095$ &$1.19095$ \\
\hline
$\mathcal{F}_5^{(4,6)}$      &$1.34255$ &$1.26048$ &$0.95234$ \\
\hline
$\mathcal{F}_6^{(4,5)}$      &$1.34255$ &$0.95234$ &$1.26048$ \\
\hline
$\mathcal{F}_6^{(4,6)}$      &$1.35938$ &$1.01481$ &$1.01481$ \\
\hline
$\mathcal{F}_7^{(3,7)}$      &$1.26829$ &$1.49903$ &$1.49903$ \\
\hline
$\mathcal{F}_7^{(3,8)}$      &$1.28228$ &$1.53709$ &$1.22995$ \\
\hline
\end{tabular}
\end{center}

Notice that if $u=w$, than the optimal configuration belongs to congruent covering ($h_1=h_2$). 

\subsubsection{On non-extendable congruent and non-congruent hyperball coverings to parameters $\{u,3,7\}$, $(6<u<7,u\in\mathbb{R})$}

We can investigate the coverings related to orthoschemes $\mathcal{F}_u^{(3,7)}$, $(6<u<7,u\in\mathbb{R})$, and here we consider the coverings in general, i.e. both congruent and non-congruent cases.

In this case, the images of the orthoschemes under reflections on its facets 
don't fill the hyperbolic space, but it provides a thinner, local covering, 
which can not be extended to the entrire $\mathbb{H}^3$. Here the previous computations are the same, 
and density function is two-variable, whose minimum 
can be determined numerically after accurate analysis.

\begin{theorem}
The non-congruent hypersphere coverings in $\mathcal{C}_u^{(3,7)}$, $(6<u<7,u\in\mathbb{R})$ attain their minimum density
at $u\approx6.45953$ with density $\approx1.26454$ where the heights of hypersperes are $h_1\approx1.50377$, and $h_2\approx1.26423$ and their common point lies on $A_1A_2$ edge. 
\end{theorem}
\begin{rmrk}
\begin{enumerate}
\item Notice, that the parameter $u$, where the minimum is attained $\approx 6.45953$, is very close to the corresponding parameter in \cite{ESz}, where it was $\approx6.45962$.
\item The above covering density is smaller than the $\approx 1.280$ density belonging to 
{\sc Fejes T\' oth L\' aszl\' o, B\" or\" oczky K\' aroly}, but this hyperball covering can not be extended to the entire space.
\end{enumerate}
\end{rmrk}
\begin{figure}[h!]
\centering
\includegraphics[width=6.5cm]{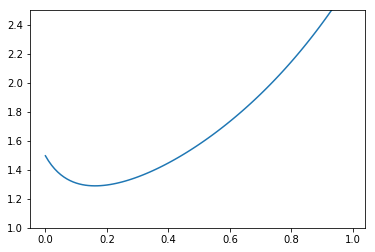} 
\includegraphics[width=6.5cm]{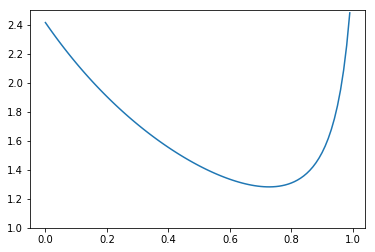}

a) \hspace{5cm} b)
\caption{a)~$\delta(\mathcal{C}_3^{(7,3)})(t)$, when $T$ lies on $QA_2$ edge ~b)~$\delta(\mathcal{C}_7^{(3,8)})(t)$, when $T$ lies on $A_1A_2$ edge~}
\end{figure}
\subsection{Congruent coverings}
In this subsection, we consider the congruent coverings for the possible cases described in Subsection 3.3.

In this case, we are looking for the $T$ point on $A_1A_2$, $CA_1$, $QA_2$ edges, which is equal distance from $HLC$ and $QEJ$. 
Investigating the above distance function, this point doesn't exists on $CA_1$ and $QA_2$, and on the third edge, 
we can find it by solving the equation $d_{HLC}^{A_1A_2}(t)=d_{QEJ}^{A_1A_2}(t)$, see Fig.~4 for visualization. 
In view of $T$, we can determine the datas of the covering as above. Here we listed just finitely many types of orthoschemes as well, but as above, 
the further types can't provide smaller densities as we obtain it after careful analysis of the density function.

\begin{center}
\begin{tabular}{|c|c|c|c|}
\hline
\multicolumn{1}{|c}{Type of orthoscheme} & \multicolumn{1}{|c|}{$\delta_{min}$} & \multicolumn{1}{|c|}{$h_1$} & \multicolumn{1}{|c|}{$h_2$}\\
\hline
$\mathcal{F}_3^{(7,3)}$      &$1.38712$ &$1.36405$ &$1.36405$ \\
\hline
$\mathcal{F}_3^{(8,3)}$      &$1.45345$ &$1.15039$ &$1.15039$ \\
\hline
$\mathcal{F}_4^{(5,4)}$      &$1.36411$ &$1.16974$ &$1.16974$ \\
\hline
$\mathcal{F}_4^{(6,4)}$      &$1.45714$ &$0.99583$ &$0.99583$ \\
\hline
$\mathcal{F}_5^{(4,5)}$      &$1.31751$ &$1.19095$ &$1.19095$ \\
\hline
$\mathcal{F}_5^{(4,6)}$      &$1.45345$ &$1.13375$ &$1.13375$ \\
\hline
$\mathcal{F}_6^{(4,5)}$      &$1.45345$ &$1.13375$ &$1.13375$ \\
\hline
$\mathcal{F}_6^{(4,6)}$      &$1.35938$ &$1.01481$ &$1.01481$ \\
\hline
$\mathcal{F}_7^{(3,7)}$      &$1.26829$ &$1.49903$ &$1.49903$ \\
\hline
$\mathcal{F}_7^{(3,8)}$      &$1.36586$ &$1.39916$ &$1.39916$ \\
\hline
\end{tabular}
\end{center}
\begin{figure}[h!]
\centering
\includegraphics[width=11cm]{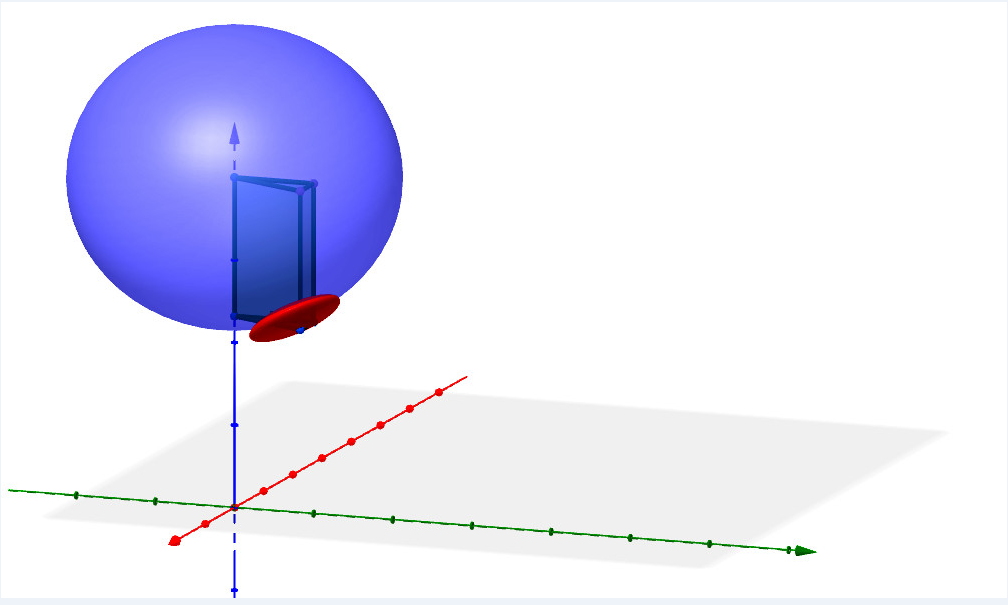} 

\caption{Visualization of $\mathcal{C}_7^{(3,7)}$ with density $\approx 1.26829$}
\end{figure}
We summarized our results in the following theorem.

\begin{theorem}
In hyperbolic space $\mathbb{H}^3$, between the congruent and non-congruent hyperball coverings generated by doubly truncated Coxeter orthoschemes, 
the $\mathcal{C}_7^{(3,7)}$  congruent hyperball configuration provides the thinnest covering with density 
$\approx 1.26829$, which is the so far known smallest ball covering density in $\mathbb{H}^3$.
\end{theorem}

\section{Hypercycle covering in hyperbolic plane}
In the hyperbolic plane ($n=2$) {\sc I. Vermes} proved, that the lower bound for congruent hypercycle covering density 
is $\frac{\sqrt{12}}{\pi}$ \cite{VI2}. However, there are no results related to thinnest hypercycle coverings with non-congruent hypercycles.
The investigation of the density of non-congruent hypersphere coverings is generally not easy. 

Here we will prove, using the results of our paper \cite{ESz}, that the theoretic lower bound $\frac{\sqrt{12}}{\pi}$ for congruent hypercycle coverings can be arbitrary approximated with non-congruent hypercycle coverings 
related to doubly truncated orthoschemes. 
\begin{theorem}
Let $A(1,0,a)$ and $B(1,b,0)$ outer points related to the Cayley-Klein-Beltrami circle model (see Fig.~5), and the $AB$ line intersects the model circle.
Let the base lines of two hypercycles be $OE$ and $FC$, and the corresponding hypercycles through the midpoint $J$ of segment $CD$ generate a covering configuration
$\mathcal{C}_{a,b} \left( \frac{1}{2}\right)$ in truncated orthosceme $FCDEO$. Then  
\begin{equation*}
\lim_{(a,b)\rightarrow (1,\infty)} \delta \left( \mathcal{C}_{a,b}\left(\frac{1}{2}\right)\right) = \frac{\sqrt{12}}{\pi}
\end{equation*} 
\end{theorem}
{\bf Proof:}

\begin{figure}[h!]
\centering
\includegraphics[width=12cm]{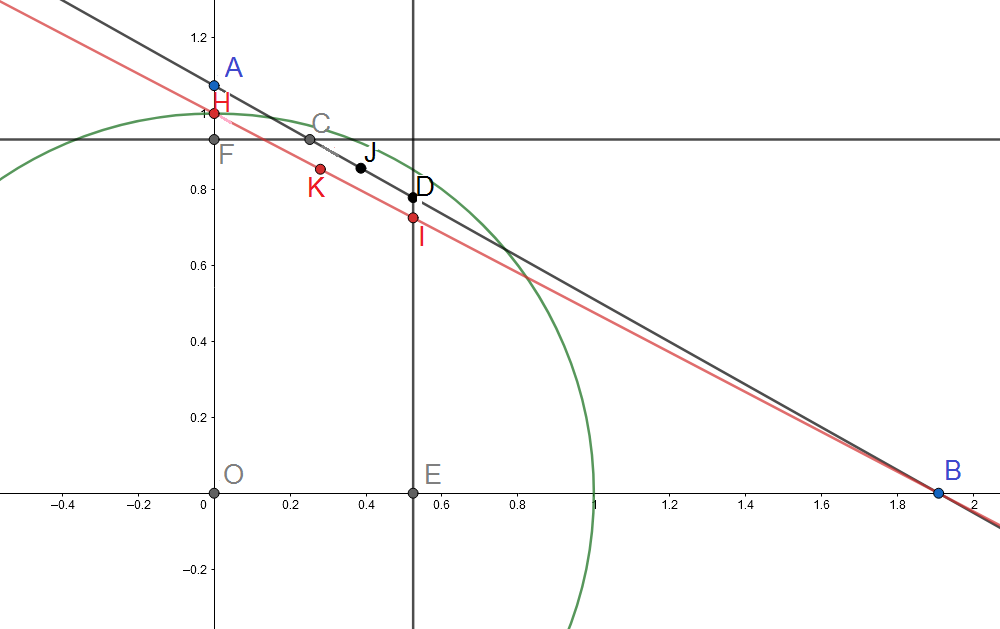} 

\caption~$FCDEO$ doubly truncated Coxeter orthoscheme~
\end{figure}

Here let the base hyperlines of the two hypercycles $OE$ and $FC$, and both hypercycles pass through the point $J$ that is the midpoint 
of $CD$. This yields obviously a covering, denote it with $\mathcal{C}_{a,b}\left(\frac{1}{2}\right)$. 
If $a\rightarrow 1$ then $A\rightarrow H$, $CD\rightarrow HI$, and the $K$ point on $HI$ arises as the limit of the $J$ midpoints of $CD$. 
Hence the orthoscheme tends to $HIEO$ that is a simple assymptotic, simple truncated orthoscheme. 
The hypercycle, whose base line is $CF$, tends to a horocycle centred at $H$. 
The hypercycle is produced as the images of a point under reflections on lines perpendicular to a certain line $CF$, 
i.e. lines passing through the polar point ($A$) of a certain line $CF$. 
The horocycle is produced as images of a point under reflections on paralell lines, i.e. lines passing through a point at the infinity ($H$). 
It means, that if $A$ tends to $H$, than the hypercycle tends to a horocycle.

So if $a\rightarrow 1$, we get a hyp-hor covering, investigated in \cite{ESz}. We recall Theorem 3.2 of \cite{ESz} using the denotation of the present paper:
\begin{theorem}[\cite{ESz}]
Let $B(1,b,0)$ outer point of the simple truncated orthoscheme, and $\mathcal{C}_{b}\left(\frac{1}{2}\right)$ denote 
the hyp-hor covering of simple truncated orthoscheme $HIEO$, with cycles passing through the point $K$. Then 
\begin{equation*}
\lim_{b \rightarrow \infty} \delta \left( \mathcal{C}_b^1\left(\frac{1}{2}\right) \right) =\frac{\sqrt{12}}{\pi}
\end{equation*}
and $\delta \left(\mathcal{C}_b^1\left(\frac{1}{2}\right)\right)>\frac{\sqrt{12}}{\pi}$ for parameter $b > 1$.
\end{theorem}

And now according to the results of this theorem, we otain the following:

\begin{equation*}
\lim_{(a,b)\rightarrow (1,\infty)} \mathcal{C}_{a,b}\left(\frac{1}{2}\right) = \lim_{b\rightarrow \infty} \lim_{a\rightarrow 1} \mathcal{C}_{a,b}\left(\frac{1}{2}\right) 
= \lim_{b\rightarrow \infty} \mathcal{C}_{b}^1 \left(\frac{1}{2}\right) = \frac{\sqrt{12}}{\pi},
\end{equation*} 

that means, that density $\frac{\sqrt{12}}{\pi}$ can be arbitrary approximated with hypercycle coverings related to doubly truncated orthoschemes. \ \ $\square$ 
\begin{rmrk}
The general investigation of planar non-congruent hypersphere packings and coverings related to the doubly truncated orthoschemes is discussed in a forthcoming paper. 
\end{rmrk}


\noindent
\footnotesize{Budapest University of Technology and Economics Institute of Mathematics, \\
Department of Geometry, \\
H-1521 Budapest, Hungary. \\
E-mail:~szirmai@math.bme.hu \\
http://www.math.bme.hu/ $^\sim$szirmai}

\end{document}